\documentclass{article}
\usepackage{amssymb, amsmath}
\usepackage{graphics}
\input xy
\xyoption{all}

\begin{document}

\title{Howe duality for the metaplectic group acting on symplectic spinor valued forms}
\author{Svatopluk Kr\'ysl \footnote{{\it E-mail address}: Svatopluk.Krysl@mff.cuni.cz}\\ {\it \small  Charles University, Sokolovsk\'a 83, Praha 8 - Karl\'in, Czech Republic} \\{\it \small and
Humboldt-Universit\"{a}t zu Berlin, Unter den Linden 6, Berlin, Germany.}
\thanks{I am very grateful to Roger Howe for an explanation of a general framework of the studied type of duality.   The author of this article was supported by the grant GA\v{C}R 306-33/80397 of the Grant Agency of the Czech Republic.
  Supported also by the SFB 647 "Space - Time - Matter" of the DFG and the Volkswagen Foundation.}}
\maketitle

\noindent\centerline{\large\bf Abstract} Let $\mathbb{S}$ denote the oscillatory module over the complex symplectic Lie algebra $\mathfrak{g}= \mathfrak{sp}(\mathbb{V}^{\mathbb{C}},\omega).$ Consider the $\mathfrak{g}$-module $\mathbb{W}=\bigwedge^{\bullet}(\mathbb{V}^*)^{\mathbb{C}}\otimes \mathbb{S}$ of exterior forms with values in the 
oscillatory module. We prove that the associative  algebra $\hbox{End}_{\mathfrak{g}}(\mathbb{W})$ is generated by the image of a certain representation of the  ortho-symplectic Lie super algebra $\mathfrak{osp}(1|2)$ and two distinguished projection operators.  
The space $\bigwedge^{\bullet}(\mathbb{V}^*)^{\mathbb{C}}\otimes \mathbb{S}$ is decomposed with respect to the joint action of  
$\mathfrak{g}$ and $\mathfrak{osp}(1|2).$ This establishes a Howe type duality for $\mathfrak{sp}(\mathbb{V}^{\mathbb{C}},\omega)$ acting on $\mathbb{W}.$

{\it Math. Subj. Class.:} 22E46, 22E47, 22E45, 81S10


\section{Introduction}

 Let $(\mathbb{V},\omega)$ be a finite dimensional real symplectic  vector space and let
   $\vartheta: Mp(\mathbb{V},\omega)\to Sp(\mathbb{V},\omega)$ denote the $2:1$ covering of the symplectic group $G = Sp(\mathbb{V},\omega)$ by the metaplectic group $\tilde{G} = Mp(\mathbb{V},\omega).$ Let us denote the the maximal compact subgroup of $\tilde{G}$ by $K,$ and the complexification of 
the Lie algebra of $G$  by $\mathfrak{g}.$  The complexification of the Lie algebra of the metaplectic group $\tilde{G}$ is isomorphic to $\mathfrak{g}$
 and thus, we shall denote it by the same letter.

 There exists a distinguished unitary  and faithful representation of the metaplectic group $\tilde{G}$ -- the so called  
Segal-Shale-Weil representation.\footnote{Let us note that also the names oscillatory,  metaplectic or symplectic spinor representation are used in the literature.} We shall denote it by $\bf S$ in this text. Its elements  are usually called symplectic spinors. For a justification
 of this name, see  Kostant \cite{Kostant}. 
Let $\mathbb{S}$ denote the underlying 
Harish-Chandra $(\mathfrak{g},K)$-module of 
$\bf S.$ When we   consider this $(\mathfrak{g},K)$-module with its $\mathfrak{g}$-module structure only, we call it the {\it oscillatory module}. 
Notice that it is known that $\mathbb{S}$ splits into two irreducible submodules $\mathbb{S} \simeq \mathbb{S}_+\oplus \mathbb{S}_-.$

 Let us set $\mathbb{W}=\bigwedge^{\bullet}(\mathbb{V}^*)^{\mathbb{C}}\otimes \mathbb{S}$ and consider the 
appropriate tensor product representation of $\mathfrak{g}$ on $\mathbb{W}$  denoting it by $\rho$. 
 In this paper, we first decompose the module $\mathbb{W}$ into irreducible 
$\mathfrak{g}$-modules. 

We shall also find generators of the 
commutant algebra $$\mbox{End}_{\mathfrak{g}}(\mathbb{W})=\{T\in \mbox{End}(\mathbb{W}) \, | \,  T \rho(X)=\rho(X)T \mbox{ for all } X\in \mathfrak{g}\}$$ 
of the symplectic Lie algebra $\mathfrak{g}$ acting on $\mathbb{W}.$
 Let $p_{\pm}: \mathbb{S} \to \mathbb{S}_{\pm}$ be the unique $\mathfrak{g}$-equivariant projections. 
These projections induce projection operators acting on the whole $\mathbb{W}$ in an obvious way. We keep denoting them by $p_{\pm}.$
We will introduce a representation $\sigma: \mathfrak{osp}(1|2)\to \mbox{End}(\mathbb{W})$ 
of the ortho-symplectic  super Lie  algebra $\mathfrak{osp}(1|2)$ on $\mathbb{W}$ and prove that 
the image of $\sigma$ together with $p_+$ and $p_-$ generate the commutant. 

At the end we decompose the $(\mathfrak{g}\times \mathfrak{osp}(1|2))$-module $\mathbb{W}$ 
into a direct sum $\bigoplus_{j=0}^l (\mathbb{E}_{jj}^-\oplus \mathbb{E}_{jj}^+) \otimes \mathbb{G}^{j},$ 
where $\mathbb{E}_{jj}^{\pm}$ is a certain irreducible infinite dimensional highest weight $\mathfrak{g}$-module
and $\mathbb{G}^j$ is an irreducible finite dimensional $\mathfrak{osp}(1|2)$-module. This establishes a Howe type duality in this case. 
   
 The basic tool used to obtain the results was the decomposition of the $\mathfrak{g}$-module $\mathbb{W}$ into irreducible summands. To obtain this decomposition, 
a theorem of Britten, Hooper, Lemire \cite{BHL} on the decomposition of the tensor product of an irreducible finite dimensional 
$\mathfrak{sp}(\mathbb{V}^{\mathbb{C}},\omega)$-module and the oscillatory module $\mathbb{S}$ was crucial.


To give a 'pure' context of our research, let us remark that the Howe dualities are generalizations of the classical results of  Schur and Weyl. Whereas Schur studied the case of $GL(\mathbb{V})$ acting on the $k$-fold product $\bigotimes^k \mathbb{V},$ Weyl considered the $SO(\mathbb{V})$-module $\bigotimes^k \mathbb{V},$ $k\in \mathbb{N}.$
See Howe \cite{Howe} for a historical treatment of the cases studied by Schur and Weyl and for some of their generalizations
(examples of what one calls  Howe type dualities nowadays). In Howe \cite{Howe}, one can also find  several applications of
these dualities and a classical version of our $2:1$ duality. (The term $2:1$ will be explained below.)

  Let us also remark that a result related to the one presented here, is the result of Slupinski in \cite{Slupinski}. In his paper, Slupinski considers the case of (orthogonal) spinor valued exterior forms as a module over the appropriate spin group. Roughly speaking, he proves that $\mathfrak{sl}(2,\mathbb{C})$ is the 
Howe dual partner to the spin group. One may rephrase this fact by saying that the situation studied in Slupinski is super symmetric to the one we are studying here.
 
 The motivation for our study of the Howe duality for exterior forms with values in the oscillatory module (i.e., for $\mathbb{W}$) comes from differential geometry and mathematical physics.  The duality we have proved here, makes it possible to do many computations in the so called symplectic spin geometry without using coordinates.  In this way, it sheds light into these structures. See, e.g., Habermann, Habermann \cite{HH} or Kr\'ysl \cite{SK2} for more details and examples. 
For applications of symplectic spinors in mathematical physics we refer an interested reader to Shale \cite{Shale}
 where they were used in a quantization procedure for Klein-Gordon fields, and to Kostant \cite{Kostant}
 who used them in geometric quantization of quantum mechanics.

 In the second section of the paper, we introduce basic notation,  summarize  basic facts on the oscillatory module and 
derive the mentioned decomposition of $\mathbb{W} = \bigwedge^{\bullet}(\mathbb{V}^{*})^{\mathbb{C}} \otimes \mathbb{S}$ 
into irreducible $\mathfrak{g}$-submodules (Theorem 3). The generators of $\mbox{End}_{\mathfrak{g}}(\mathbb{W})$  are given in the 
third section (Theorem 8). In that section, 
the representation $\sigma: \mathfrak{osp}(1|2)\to \mbox{End}(\mathbb{W})$ is introduced and the fact,
 that it is a representation is proved (Theorem 9).
 In the fourth section, the space $\mathbb{W}$ is  explicitly decomposed into irreducible
 submodules with respect to the joint action of $\mathfrak{g}$ and $\mathfrak{osp}(1|2)$  (Theorem 13).
 
\section{Decomposition of $\bigwedge^{\bullet}(\mathbb{V}^*)^{\mathbb{C}}\otimes \mathbb{S}$}

Suppose that $\mathfrak{g}$ is a simple complex Lie algebra. Choose a Cartan subalgebra $\mathfrak{h}$ of $\mathfrak{g}$ and a set of positive roots $\Phi^+.$ We denote the complex irreducible highest weight $\mathfrak{g}$-module with a highest weight $\mu \in \mathfrak{h}^*$ by $L(\mu).$ If $\mu$ happens to be dominant and integral with respect to the choice
$(\mathfrak{h}, \Phi^+),$ we denote $L(\mu)$ by $F(\mu)$ emphasizing the fact that the module $L(\mu)$ is finite dimensional.  For an integral dominant weight $\mu$ with respect to $(\mathfrak{h},\Phi^+),$ we denote the set of all weights of the irreducible representation $F(\mu)$ by $\Pi(\mu).$

Now, let us restrict our attention to the studied symplectic case.
Consider a $2l$ dimensional real symplectic vector space
$(\mathbb{V},\omega).$ Let $\mathbb{V} = \mathbb{L}\oplus \mathbb{L}'$ be a direct sum decomposition of the vector space $\mathbb{V}$ into two Lagrangian subspaces $\mathbb{L}$ and $\mathbb{L}'$ of $(\mathbb{V},\omega).$ Let $\{e_i\}_{i=1}^{2l}$ be an adapted  symplectic basis of $(\mathbb{V},\omega),$ i.e., $\{e_i\}_{i=1}^{2l}$ is a symplectic basis of $(\mathbb{V},\omega)$ and $\{e_i\}_{i=1}^{l}\subseteq \mathbb{L}$ and $\{e_{i}\}_{i=l+1}^{2l}\subseteq \mathbb{L}'.$ Because the notion of a symplectic basis is not unique, let us fix one for a later use.
We call a basis $\{e_i\}_{i=1}^{2l}$ of $\mathbb{V}$ symplectic basis of $(\mathbb{V},\omega)$ if for $\omega_{ij}=\omega(e_i,e_j),$ we have 
\begin{itemize}
\item[] $\omega_{ij}=1$ if an only if $i\leq l$ and $j=i+l,$
\item[] $\omega_{ij}=-1$ if and only if $i>l$ and $j=i-l$ and
\item[] $\omega_{ij}=0$ in other cases.
\end{itemize}
For $i,j =1,\ldots, 2l,$ let us define numbers $\omega^{ij}$ by $\sum_{k=1}^{2l}\omega_{ik}\omega^{jk}=\delta_i^j.$ 
The basis of $\mathbb{V}^*$ dual to the basis $\{e_i\}_{i=1}^{2l}$
will be denoted by $\{\epsilon^i\}_{i=1}^{2l}.$ 


Let us denote the symplectic group $Sp(\mathbb{V},\omega)$ by $G.$ In this text, the metaplectic group $Mp(\mathbb{V},\omega)$ is denoted by $\tilde{G}$ and the 
appropriate two-fold covering map $\tilde{G} \to G$ by $\vartheta.$
For some technical reasons, let us suppose that $l\geq 2,$  and consider this choice holds throughout the rest of this article.  We shall denote the complexification of the Lie algebra of $(\mathbb{V},\omega),$ i.e., the Lie algebra  $\mathfrak{sp}(\mathbb{V}^{\mathbb{C}},\omega)$ by $\mathfrak{g}.$
(The complexified symplectic form on $\mathbb{V}^{\mathbb{C}}$ 
will still be denoted by $\omega$.) If we choose a Cartan subalgebra
$\mathfrak{h}\subseteq \mathfrak{g}$ and a
system of positive roots $\Phi^+ \subseteq
\mathfrak{h}^*,$ then the set of fundamental weights
$ \{ \varpi_i \}_{i=1}^l$ is uniquely determined. Further this set
determines a basis $\{\epsilon_i\}_{i=1}^l$ of
$\mathfrak{h}^*$ given by
$\varpi_i=\sum_{j=1}^i\epsilon_j$ for $i=1,\ldots, l.$\footnote{Let
us mention that we are distinguishing the sets
$\{\epsilon^i\}_{i=1}^{l} \subseteq \mathbb{V}^*$ and
$\{\epsilon_i\}_{i=1}^l \subseteq \mathfrak{h}^{*}$ from each other,
although we could have identified them by a choice of the first
inclusion from the chain of inclusions $\mathfrak{h}\subseteq \mathfrak{g} \subseteq (\mathbb{V}^{*}\otimes \mathbb{V})^{\mathbb{C}}$ and 
a set of positive roots.}
For $\mu=\sum_{i=1}^l\mu_i\epsilon_i,$ we shall often denote $L(\mu)$ by $L(\mu_1,\ldots,\mu_l)$ simply.


The Segal-Shale-Weil representation is a faithful unitary representation of the metaplectic group
$Mp(\mathbb{V},\omega)$ on the complex vector  space $L^2(\mathbb{L})$ of complex valued square Lebesgue integrable functions defined on $\mathbb{L}.$ Because we would like to omit problems caused by dealing with unbounded operators, we shall consider the underlying Harish-Chandra $(\mathfrak{g},K)$-module of the Segal-Shale-Weil representation. 
This module considered as a $\mathfrak{g}$-module only will be called the {\it oscillatory module} and the appropriate representation
$$L: \mathfrak{g} \to \mbox{End}(\mathbb{S})$$ will be called the {\it oscillatory representation} (see Howe \cite{Howe}).
One can prove that  $\mathbb{S}_+ \simeq 
L(\mu^0)$  and $\mathbb{S}_-\simeq L(\mu^1)$ where $\mu^0= -\frac{1}{2}\varpi_l$ and $\mu^1=\varpi_{l-1}-\frac{3}{2}\varpi_l.$  Let us notice that as a vector space 
$\mathbb{S}_+$ is isomorphic to the space of even polynomials in $\mathbb{C}[z^1,\ldots, z^l]$ and $\mathbb{S}_-$ to the space of the odd ones.
For more information on the Segal-Shale-Weil representation, see
Weil \cite{Weil} and Kashiwara, Vergne \cite{KV}. For information on the oscillatory module, see Britten, Hooper, Lemire \cite{BHL}. 

 
%
 
In order to derive the studied Howe type duality, we shall need the symplectic Clifford multiplication $\mathbb{V}^{\mathbb{C}}\times \mathbb{S} \to \mathbb{S}$, which enables us to multiply symplectic spinors by  element from $\mathbb{V}^{\mathbb{C}}.$ It is given by the following prescription 
\begin{eqnarray}
(e_i.f)(x)=\frac{\partial f}{\partial x^i}(x), \, \mbox{ } (e_{i+l}.f)(x)=\imath x^i f(x),\, i=1,\ldots, l \label{sc}
\end{eqnarray} where $x = \sum_{i=1}^lx^ie_i \in \mathbb{L}, f\in \mathbb{S}$ and it is extended linearly to the whole $\mathbb{V}^{\mathbb{C}}.$
(The symplectic Clifford multiplication is basically the canonical  quantization prescription.) 


Now, for $i = 0, 1$ and a dominant integral weight $\mu=\sum_{j=1}^l\mu_j \varpi_j \in \mathfrak{h}^{*},$ let us introduce two sets $T^i_{\mu}$.  A weight 
$\mu' \in \mathfrak{h}^*$ is an element of $T^i_{\mu}$ if and only if the numbers $d_j,$ $j=0,\ldots, l,$ defined by 
$\mu'-\mu = \sum_{j=1}^ld_j\epsilon_j,$ satisfy the following  conditions
\begin{itemize}	
\item[1)] $d_j + \delta_{l,j}\delta_{1,i} \in \mathbb{N}_0$ for $j=0,\ldots,l$
\item[2)] $0\leq d_j \leq \mu_j$ for $j=0,\ldots, l - 1$ and  $0 \leq d_l +\delta_{1,i} \leq 2\mu_l +1$ 
\item[3)] $\sum_{j=1}^l d_j$ is even. 
\end{itemize}
 

In what follows, we shall  need a result on the decomposition of the tensor product of a finite dimensional $\mathfrak{g}$-module with one of the modules $L(\mu^i),$ $i=0,1,$ into irreducible submodules.  This result was published in Britten, Hooper, Lemire \cite{BHL}. Because there is a misprint in the cited article, let us  write the corrected result explicitly here.

{\bf Theorem 1:}
For $i=0,1$ and an integral dominant weight $\mu$ with respect to $(\mathfrak{h},\Phi^+),$ we have
$$F(\mu) \otimes L(\mu^i) \, \simeq  \bigoplus_{\kappa \in T_{\mu}^i \cap \Pi(\mu)}L(\mu_i + \kappa).$$

{\it Proof.} See Britten, Hooper, Lemire \cite{BHL}. The mentioned misprint occurs in the definition of $T_{\mu}^1$ in that paper (private communication with F. W. Lemire). $\Box$

For convenience let us introduce a function $\mbox{sgn}:\{+,-\}\to \{0,1\}$ given by the prescription
$\mbox{sgn}(+)=0$ and $\mbox{sgn}(-)=1.$  Following modules will be used in the decomposition theorem
$$\mathbb{E}_{ij}^{\pm}:=L(\underbrace{\frac{1}{2},\cdots,\frac{1}{2}}_{j},\underbrace{-\frac{1}{2},\cdots,
-\frac{1}{2}}_{l-j-1},-1+\frac{1}{2}(-1)^{i+j+\mbox{sgn}(\pm)})$$ where $i=0,\ldots,l-1,$   $j=0,\ldots, i$ and $i=l,$   $j=0,\ldots, l-1.$ For $i=j=l,$
we set $\mathbb{E}^+_{ll}:=L(\frac{1}{2}, \cdots, \frac{1}{2})$ and $\mathbb{E}^-_{ll}:=L(\frac{1}{2}, \cdots, \frac{1}{2}, -\frac{5}{2}).$
For $i=l+1,\ldots,2l$ and $j=0,\ldots, 2l-i,$ we set $\mathbb{E}_{ij}^{\pm}:=\mathbb{E}_{(2l-i)j}^{\pm}.$ 
In order to write the results  as short as possible, let us introduce the following notation.
For $i=0,\ldots, l,$ let us define $m_i=i$ and for $i=l+1,\ldots, 2l,$   $m_i=2l-i.$ With these conventions, let  us set 
$$\Xi=\{(i,j)|i=0,\ldots, 2l; j= 0,\ldots, m_i\} \, \mbox{  and further  } 
 \mathbb{E}_{ij}^{\pm} = 0 \mbox{ for } (i,j)\in \mathbb{Z}^2 \setminus \Xi.$$

Now, let us derive the next

{\bf Lemma 2:} For $r=1,\ldots, l,$
we have $$\Pi(\varpi_r) \supseteq \{\sum_{s=1}^r\pm\epsilon^{i_s}|1\leq i_1 < \ldots < i_r  \leq l\}.$$

{\it Proof.} According to the Corollary  5.1.11. pp. 237 and Theorem 5.1.8. (3) pp. 236 in Goodman, Wallach \cite{GW}, the
$\mathfrak{g}$-module $F(\varpi_r)$ is isomorphic to the linear span
of isotropic $r$-vectors in $\mathbb{V}^{\mathbb{C}}$, i.e., multi-vectors $v=f_1\wedge\ldots\wedge f_r,$ where $\omega(f_i,f_j)=0$ for $i,j=1,\ldots, l,$ $r=1,\ldots, 2l.$ Second, it is easy to realize that one can choose a Cartan subalgebra $\mathfrak{h}$ of $\mathfrak{g}$ in a way that the following is true. For $i=1,\ldots, l,$ the basis vector $e_i \in \mathbb{V}^{\mathbb{C}}$ is a weight vector of the weight $\epsilon_i$ and the  vector $e_{i+l}$ is a  weight vector of the weight $-\epsilon_i,$ both for the defining representation $\mathfrak{g}$ on $\mathbb{V}^{\mathbb{C}}.$ Using this, the result follows easily.
$\Box$

Now, we  introduce the module $\mathbb{W}$ we shall be mostly concern with in this paper. As a vector space
$$\mathbb{W}=\bigwedge^{\bullet}(\mathbb{V}^{*})^{\mathbb{C}}\otimes \mathbb{S}.$$ 
The representation $ \rho: \mathfrak{g} \to \mbox{End}(\mathbb{W})$ of $\mathfrak{g}$ on ${\mathbb W}$ we are interested in, is defined by the prescription $$\rho(X)(\alpha \otimes s):= \mbox{ad}(\mbox{X})^{* \wedge r}(\alpha) \otimes L(X)s$$ for
 $X\in \mathfrak{g},$ $\alpha\in \bigwedge^r (\mathbb{V}^*)^{\mathbb{C}},$ $r=0,\ldots, 2l$ and $s\in \mathbb{S}.$ The module $\mathbb{W}$ will be   called the
module of (exterior) forms with values in symplectic spinors.

Now, we can state the decomposition theorem. The proof of this theorem is based on a direct application of Theorem 1 and Lemma 2.

{\bf Theorem 3:} For $i=0,\ldots,2l,$ the following decomposition into irreducible $\mathfrak{g}$-modules
$$\bigwedge^i (\mathbb{V}^{*})^{\mathbb{C}}\otimes \mathbb{S}_{\pm} \simeq \bigoplus_{j, (i,j)\in\Xi} \mathbb{E}_{ij}^{\pm} \quad \mbox{ holds.}$$

{\it Proof.}  
Using Theorem 5.1.8. pp. 236 and Corollary 5.1.9. pp. 237 in Goodman, Wallach \cite{GW}, we have for
$i=2k,$ $k\in \mathbb{N}_0$
\begin{eqnarray}
\bigwedge^i (\mathbb{V}^*)^{\mathbb{C}}\otimes \mathbb{S}_{\pm}=(F(\varpi_0)\oplus F(\varpi_2)\oplus \ldots \oplus F(\varpi_i))\otimes \mathbb{S}_{\pm}. \label{1}
\end{eqnarray}
Using the cited theorems in Goodman, Wallach \cite{GW}  again, we have have for $i= 2k+1,$ $k \in \mathbb{N}_0$
\begin{eqnarray}
\bigwedge^i (\mathbb{V}^*)^{\mathbb{C}}\otimes \mathbb{S}_{\pm}=(F(\varpi_1)\oplus F(\varpi_3)\oplus \ldots \oplus F(\varpi_i))\otimes \mathbb{S}_{\pm}. \label{2}
\end{eqnarray}

We shall consider the mentioned tensor products  for $i=0,\ldots,l$ only, because the result for  $i=l+1,\ldots,2l$ follows from that one for $i=0,\ldots,l$ immediately
due to the  $\mathfrak{g}$-isomorphism $\bigwedge^i (\mathbb{V}^*)^{\mathbb{C}}\otimes \mathbb{S}_{\pm} \simeq \bigwedge^{2l-i} (\mathbb{V}^*)^{\mathbb{C}}\otimes \mathbb{S}_{\pm}$ and the definition of $\mathbb{E}_{ij}^{\pm}$ for
$i=l+1,\ldots,2l$ and $j=0,\ldots, m_i.$

Let us consider the tensor products by $\mathbb{S}_+$ and $\mathbb{S}_-$ separately.

\begin{itemize}
\item[1)] First, let us consider the tensor product $\bigwedge^i (\mathbb{V}^*)^{\mathbb{C}}\otimes \mathbb{S}_+$
for $i=0,\ldots,l.$

Using Lemma 2 and Theorem 1, we can easily compute that for $j=1,\ldots, l,$
$T_{\varpi_j}^0=\{\epsilon_1+\ldots+\epsilon_j, \epsilon_1 + \ldots +\epsilon_{j-1}- \epsilon_l\} \subseteq \Pi(\varpi_j)$ and thus,
$F(\varpi_j)\otimes \mathbb{S}_+=L(\underbrace{\frac{1}{2},\ldots,\frac{1}{2}}_j,\underbrace{-\frac{1}{2},\ldots,-\frac{1}{2}}_{l-j})\oplus
L(\underbrace{\frac{1}{2},\ldots,\frac{1}{2}}_{j-1},\underbrace{-\frac{1}{2},\ldots,-\frac{1}{2}}_{l-j},-\frac{3}{2})$ using the relation $\omega_j=\sum_{i=1}^j\epsilon_i.$
Summing up these terms according to (\ref{1}) and (\ref{2}), we obtain the statement of the theorem for both of the cases $i$  is odd and $i$ is even.

\item[2)]  Now, let us consider the tensor product $\bigwedge^i(\mathbb{V}^*)^{\mathbb{C}}\otimes \mathbb{S}_-$ for $i=0,\ldots,l.$

 Using Lemma 2, we can easily compute that for $j=1,\ldots, l-1,$ we have
$T_{\varpi_j}^1=\{\epsilon_1+\ldots+\epsilon_j,\epsilon_1+\ldots+\epsilon_{j-1}+\epsilon_l\}\subseteq \Pi(\omega_j)$ and $T^1_{\varpi_l}=\{\epsilon_1+\ldots + \epsilon_l, \epsilon_1+\ldots + \epsilon_{l-1}-\epsilon_l\}\subseteq \Pi(\varpi_l).$ 
Therefore using Theorem 1, we have 
$F(\omega_j)\otimes \mathbb{S}_-=L(\underbrace{\frac{1}{2},\ldots,\frac{1}{2}}_{j-1},\underbrace{-\frac{1}{2},\ldots,-\frac{1}{2}}_{l-j+1})\oplus L(\underbrace{\frac{1}{2},\ldots,\frac{1}{2}}_{j},\underbrace{-\frac{1}{2},\ldots,-\frac{1}{2}}_{l-j-1},-\frac{3}{2})$ for $j=1,\ldots, l-1.$ For $j=l,$ we obtain
$F(\varpi_l)\otimes \mathbb{S}_+=L(\frac{1}{2} \ldots \frac{1}{2} -\frac{3}{2})\oplus L(\frac{1}{2} \ldots \frac{1}{2}-\frac{5}{2})$ using Theorem 1 again.
Summing up these terms according to (\ref{1}) and (\ref{2}), we obtain the statement of the theorem for both cases
$i$ is odd and $i$ is even.
\end{itemize}
$\Box$

{\bf Remark:} Due to Theorem 3 and definitions of $\mathbb{E}_{ij}$, we know that for $i=0,\ldots, 2l$ the $\mathfrak{g}$-module $\bigwedge^{i}(\mathbb{V}^{*})^{\mathbb{C}}\otimes \mathbb{S}$ is 
multiplicity-free.

\section{The commutant $\hbox{End}_{\mathfrak{g}}(\mathbb{W})$ and its supersymmetric realization}

In this section, we shall prove that the associative algebra $\hbox{End}_{\mathfrak{g}}(\bigwedge^{\bullet}(\mathbb{V}^*)^{\mathbb{C}}\otimes \mathbb{S})$ is generated  by the below introduced four explicitly given elements (a "raising", a "lowering" operator and two projections).  

 For $r=0,\ldots, 2l$ and $\alpha \otimes s = \alpha \otimes (s_+, s_-) \in \bigwedge^r
(\mathbb{V}^*)^{\mathbb{C}}\otimes (\mathbb{S}_+ \otimes \mathbb{S}_-),$  we set
$$F^+:\bigwedge ^{r}(\mathbb{V}^*)^{\mathbb{C}}\otimes \mathbb{S} \to
\bigwedge^{r+1}(\mathbb{V}^*)^{\mathbb{C}}\otimes \mathbb{S},\, F^+(\alpha \otimes
s)=\frac{\imath}{2}\sum_{i=1}^{2l}\epsilon^i\wedge \alpha \otimes e_i.s,$$
$$F^-:\bigwedge ^{r}(\mathbb{V}^*)^{\mathbb{C}}\otimes \mathbb{S} \to \bigwedge ^{r-1}(\mathbb{V}^*)^{\mathbb{C}}\otimes \mathbb{S}, \,
F^-(\alpha \otimes s)=\frac{1}{2}\sum_{i=1}^{2l}\omega^{ij}\iota_{e_i}\alpha \otimes
e_j.s,$$ 
$$p_{\pm}:\bigwedge ^{r}(\mathbb{V}^*)^{\mathbb{C}}\otimes \mathbb{S} \to \bigwedge ^{r}(\mathbb{V}^*)^{\mathbb{C}}\otimes \mathbb{S}, \, p_{\pm}(\alpha\otimes s)=\alpha \otimes s_{\pm}$$
and extended linearly to the whole $\mathbb{W}.$ Next, consider the operator $H$ defined by the formula
$$H=2 (F^+ F^- + F^-F^+).$$
 
 The values of the operator
$H$ are determined in the next

{\bf Lemma 4:} Let $(\mathbb{V},\omega)$ be a 
symplectic vector space of dimension $2l.$   Then for   $r=0,\ldots, 2l,$ we have
$$H_{|\bigwedge^r(\mathbb{V}^*)^{\mathbb{C}} \otimes \mathbb{S}}=\frac{1}{2}(r-l)\hbox{Id}_{|\bigwedge^r(\mathbb{V}^*)^{\mathbb{C}}\otimes \mathbb{S}}.$$

{\it Proof.} The proof is straightforward, see Kr\'ysl \cite{SK}. $\Box$

{\bf Lemma 5:} The maps $F^{\pm}, p_{\pm}$ and $H$ are $\mathfrak{g}$-equivariant
with respect to the representation $\rho$ of $\mathfrak{g}$ on $\mathbb{W}.$

{\it Proof.} The operators $p_{\pm}$ are clearly $\mathfrak{g}$-equivariant. The $\mathfrak{g}$-equivariance of $F^{\pm}$ and $H$
can be verified straightforwardly. See Kr\'ysl \cite{SK} for a proof. 
$\Box$

{\bf Definition 1:} Let us denote the associative algebra generated by $F^{\pm}$ and $p_{\pm}$   by $\mathfrak{C}$.

Due to the previous lemma, we already know that $\mathfrak{C} \subseteq \mbox{End}_{\mathfrak{g}}(\mathbb{W}).$  Now, we shall prove that $\mathfrak{C}$ exhausts the whole 
commutant $\mbox{End}_{\mathfrak{g}}(\mathbb{W}).$ For convenience, let us define the following two subsets $\Xi_+$ and $\Xi_-$ of the set $\Xi.$ We define $\Xi_-=\Xi \setminus \{(j,2l-j)| j=l,\ldots, 2l\}$ and $\Xi_+=\Xi \setminus \{(j,j)| j=0,\ldots, l\}.$ The sets $\Xi_+$ and $\Xi_-$ can be visualized as the left hand side and the right hand part of the 'triangular shape' of $\Xi,$ respectively. 
 
{\bf Lemma 6:} For each $(i,j)\in \Xi,$ we have
\begin{eqnarray*}
F^{+}_{|\mathbb{E}_{ij}^{\pm}}&:& \mathbb{E}_{ij}^{\pm}\stackrel{\sim}{\rightarrow}
\mathbb{E}_{i + 1,j}^{\mp} \qquad \mbox{ if } (i,j) \in \Xi_{-} \mbox{ and}\\
F^{-}_{|\mathbb{E}_{ij}^{\pm}}&:& \mathbb{E}_{ij}^{\pm}\stackrel{\sim}{\rightarrow}
\mathbb{E}_{i - 1,j}^{\mp} \qquad \mbox{ if } (i,j) \in \Xi_{+}.
\end{eqnarray*}
 
{\it Proof.} First, let us prove by induction the following formulas. For $(i,j) \in \Xi$ we have
\begin{eqnarray} \label{formule}
 F^-F^+_{|\mathbb{E}_{ij}}= \left\{\begin{array}{l}
                                                       \frac{1}{4}\left(\frac{1+i-j}{2}\right)\mbox{Id}_{|\mathbb{E}_{ij}} \quad  \hbox{ if } \ i+j \hbox{ is odd}   \hbox{,} \\
                                                      \frac{1}{4}\left(\frac{i+j}{2} -l\right)\mbox{Id}_{|\mathbb{E}_{ij}}  \quad {\hbox{if} \  i+j \hbox{ is even.}
                                                       }
                                                       \end{array}
                                                       \right.
                                                       \end{eqnarray}

Let us fix an integer $j \in \{0,\ldots, l\}$ and proceed by the induction on the form degree $i.$                                                  
\begin{itemize} 
\item[I.] For $i =j $  and $\phi \in \mathbb{E}_{ii},$ let us compute $F^-F^+\phi = (\frac{1}{2} H -  F^+F^-)\phi= \frac{1}{4}(i-l)\phi- F^+F^-\phi$ due to the definition of $H$ and Lemma 4. We have $F^-_{|\mathbb{E}_{ii}}=0$ because $F^-$ is $\mathfrak{g}$-equivariant (Lemma 5), lowers the form degree by one and there is no summand isomorphic 
to $\mathbb{E}_{ii}$ in $\bigwedge^{i-1}(\mathbb{V}^{*})^{\mathbb{C}}\otimes \mathbb{S}$
 (see the Theorem 3). Summing up, we have $F^-F^+\phi=\frac{1}{4}(i-l)$ according to (\ref{formule}).
\item[II.] Now, let us suppose the statement is true for $(i,j) \in \Xi,$ $i+j$ is odd. For $(i+1,j)\in \Xi$ and $\phi \in \mathbb{E}_{i+1,j},$ let us compute $F^-F^+\phi = \frac{1}{2} H \phi - F^+F^-\phi = \frac{1}{4}(i+1-l)\phi - F^+F^-\phi$  due to the definition of $H$ and Lemma 4.
Using the induction hypothesis, we  have that $F^-F^+_{|\mathbb{E}_{ij}}=\frac{1}{4}(\frac{1+i-j}{2})\mbox{Id}_{|\mathbb{E}_{ij}}.$ Thus $F^+_{|\mathbb{E}_{ij}}$ is injective. Because $F^+$ is $\mathfrak{g}$-equivariant, 
raises the form degree by one and there is no other summand in $\bigwedge^{i+1}(\mathbb{V}^*)^{\mathbb{C}}\otimes \mathbb{S}$ isomorphic to $\mathbb{E}_{ij}$ than $\mathbb{E}_{i+1,j},$ we see that $F^+_{|\mathbb{E}_{ij}}: \mathbb{E}_{ij} \to \mathbb{E}_{i+1,j}.$ Because of the proved injectivity and the Schur lemma for complex highest weight modules (see Diximier \cite{Dixmier}), we know that $F^{+}_{|\mathbb{E}_{ij}}$ is actually an isomorphism. Thus there exists $\tilde{\phi}\in \mathbb{E}_{ij}$ such that $\phi=F^+\tilde{\phi}.$
We may compute $F^+F^-\phi=F^+F^-(F^+\tilde{\phi})=F^+(F^-F^+\tilde{\phi})=\frac{1}{4}(\frac{1+i-j}{2})F^+\tilde{\phi}=\frac{1}{4}(\frac{1+i-j}{2})\phi$ by the induction hypothesis. Substituting this relation into the already derived $F^-F^+\phi=\frac{1}{4}(i+1-l)\phi - F^+F^-\phi,$ we get $F^-F^+\phi = \frac{1}{4}(i+1-l)\phi - \frac{1}{4}(\frac{1+i-j}{2})\phi= \frac{1}{4}(\frac{i +1 + j}{2}-l)\phi$ according to the formula (\ref{formule}).

Now, let us suppose the statement is true for $(i,j) \in \Xi$ and $i+j$ is even.  For $(i+1,j)\in \Xi$ and $\phi \in \mathbb{E}_{i+1,j},$ we can compute $F^-F^+\phi=\frac{1}{2}H\phi-F^+F^-\phi = \frac{1}{4}(i+1-l)\phi-F^+F^-\phi$ due to the definition of $H$ and  Lemma 4.
By parallel lines of reasoning, we get the existence of $\tilde{\phi} \in \mathbb{E}_{ij}$ such that $\phi = F^+\tilde{\phi}.$  Using the induction hypothesis, we may write $F^+F^-\phi=F^+F^-(F^+\tilde{\phi})=F^+(F^-F^+\tilde{\phi})=\frac{1}{4}(\frac{i-j +1}{2}-l)F^+\tilde{\phi}=\frac{1}{4}(\frac{i+j}{2}-l) \phi$ by the induction hypothesis. Substituting this expression into the computation above, we get $F^-F^+\phi=\frac{1}{4}(i+1-l)\phi - \frac{1}{4}(\frac{i+j}{2}-l)\phi =\frac{1}{4}(\frac{1 + (i +1) -j}{2})\phi.$ Thus, the formula follows.
\end{itemize}

Using the formula (\ref{formule}), we see that 
$F^-F^+_{|\mathbb{E}_{ij}}$ is injective iff $i+j \neq 2l$ and $j \neq i+1$, i.e., $(i,j)\in \Xi_-,$ the second condition being empty. Especially, $F^+$ is injective for $(i,j)\in \Xi_-.$ Using the Schur lemma, we get that $F^+$ is an isomorphism. From this, we also know that $F^-$ is injective on the image of $F^+,$ i.e., it is an isomorphism for $(i,j)\in \Xi_+.$
$\Box$

{\bf Remark:} It is  easy to see that $F^-$ is zero when restricted to 
$\mathbb{E}_{ii}^{\pm},$ $i=0,\ldots,l.$  Namely, we know that $F^-$ is lowering the form degree by one, it is a $\mathfrak{g}$-equivariant map and  there is no submodule of the module of symplectic spinor valued exterior forms of form degree  $i-1$ isomorphic to $\mathbb{E}_{ii}^{\pm}$ (see Theorem 3). 
A similar discussion can be made for $F^+$ restricted to $\mathbb{E}_{im_i},$ $i=l,\ldots, 2l.$ 

For $(i,j)\in \Xi,$ let us denote the projection operators from the space $\mathbb{W}$ to the submodule
$\mathbb{E}_{ij}^{\pm}$ by $S_{ij}^{\pm},$ i.e.,
$$S_{ij}^{\pm}:\bigwedge^{\bullet}(\mathbb{V}^*)^{\mathbb{C}} \otimes \mathbb{S}_{\pm}\to \mathbb{E}_{ij}^{\pm} \subseteq \bigwedge^i (\mathbb{V}^*)^{\mathbb{C}} \otimes \mathbb{S}_{\pm}.$$ 

{\bf Lemma 7:} For each $(i,j)\in \Xi,$ the projection operators
$S_{ij}^{\pm}$ are elements of the associative algebra $\mathfrak{C} = \langle F^+, F^-, p_+, p_- \rangle.$

{\it Proof.} For each $i = 0, \ldots, 2l,$ let us  define the following operators $S_i^{\pm}:
\bigwedge^{\bullet}(\mathbb{V}^*)^{\mathbb{C}} \otimes \mathbb{S}_{\pm} \to
\bigwedge^{i}(\mathbb{V}^*)^{\mathbb{C}}\otimes \mathbb{S}_{\pm}$ by the formula
$$S_i^{\pm}:=\left( \prod_{k=0, k\neq i}^{2l}\frac{2H-k+l}{i-k}\right) p_{\pm}.$$
Using Lemma 4, we see that the image of each $S_i^{\pm}$ lies in the
prescribed space and the normalization is correct. Recall that due to its definition, $H$ can be expressed using  the operators $F^{+}$ and
$F^-$ only and thus for $i=0,\ldots, 2l,$ the projection $S_i^{\pm}$ lies in
the algebra $\mathfrak{C}$. Now let us fix an integer $i\in \{0,\ldots, 2l\}.$ We prove that for
each $j$ such that $(i,j)\in \Xi,$ the projection
$S_{ij}^{\pm}:\bigwedge^{i}(\mathbb{V}^*)^{\mathbb{C}}\otimes \mathbb{S}_{\pm} \to \mathbb{E}_{ij}^{\pm} \in \mathfrak{C}.$ 
  We proceed by
induction (on $j$).
\begin{itemize}
\item[I.] For $j=0,$ we can define $S_{i0}'':=(F^+)^{i}(F^-)^{i}.$ Using
the fact that applying the $F^-$ (or $F^{+}$) lowers (or raises) the form
degree by $1,$ we see that $S_{i0}'': \bigwedge^{i}(\mathbb{V}^*)^{\mathbb{C}}\otimes
\mathbb{S}_{\pm}\to \mathbb{E}_{i0}^{\pm}$. Using the Schur lemma 
for complex irreducible highest weight modules (see Dixmier \cite{Dixmier}), we see that there exists a complex number $\lambda_{i0}\in
\mathbb{C}$ such that ${S_{i0}''}_{|\mathbb{E}^{\pm}_{i0}}=\lambda_{i0}\mbox{Id}_{|\mathbb{E}_{i0}^{\pm}}.$ 
Due to  Lemma 6, we know that
$\lambda_{i0}\neq 0.$ Thus defining
$S_{i0}^{\pm}:=\frac{1}{\lambda_{i0}}S_{i0}''\circ S_i^{\pm},$ we get a projection from $\bigwedge^{\bullet}(\mathbb{V}^*)^{\mathbb{C}}\otimes \mathbb{S}_{\pm}$ onto
$\mathbb{E}_{i0}^{\pm}$ expressed as linear combinations of compositions of
$F^{\pm}$ and $p_{\pm}.$

\item[II.] Let us suppose that for $k=0,\ldots, j,$ the operators $S_{ik}^{\pm}$
can be written as linear combinations of compositions of the
operators $F^{\pm}$ and $p_{\pm}.$ Now, we shall use the
operators $S_{i0}^{\pm},\ldots, S_{ij}^{\pm}$ in order to define $S_{i,j+1}^{\pm}.$
Let us take an element $\xi \in \bigwedge^{i}(\mathbb{V}^*)^{\mathbb{C}}\otimes \mathbb{S}_{\pm}$
and form an element $\zeta:= S_{i,j+1}'\xi:=\xi -
\sum_{k=0}^jS_{ik}^{\pm}\xi \in \bigoplus_{k=j+1}^{m_i} \mathbb{E}_{ik}^{\pm}.$ Now, form an element
$\zeta':=S_{i,j+1}''\zeta:=(F^+)^{i-j}(F^-)^{i-j}\zeta.$ In the
same way as in the item I., we conclude that $\zeta' \in \mathbb{E}_{i,j+1}^{\pm}.$ 
Form the homomorphism $S_{i,j+1}^{'''}:=S_{i,j+1}^{''}\circ
S_{i,j+1}'.$
Using the Schur lemma in the case of $S^{'''}_{i,j+1 |\mathbb{E}_{i,j+1}^{\pm}}: \mathbb{E}_{i,j+1}^{\pm} \to \mathbb{E}_{i,j+1}^{\pm},$ we
conclude that there is a complex number $\lambda_{i,j+1} \in
\mathbb{C}$ such that $S^{'''}_{i,j+1|\mathbb{E}_{i,j+1}^{\pm}}=\lambda_{i,j+1}
\mbox{Id}_{|\mathbb{E}_{i,j+1}^{\pm}}.$  Due to Lemma 6, we know that
$\lambda_{i,j+1}\neq 0.$ Now 
$S_{i,j+1}^{\pm}:=\frac{1}{\lambda_{i,j+1}}S_{i,j+1}'''\circ S_i^{\pm}$ is the desired projection. Going through the construction back, we see that  we have used the operators $F^{\pm}$ and $p_{\pm}$ only. $\Box$
\end{itemize}
 
Now, we prove that the algebra $\mathfrak{C}$ exhausts the whole commutant $\mbox{End}_{\mathfrak{g}}(\mathbb{W}).$

{\bf Theorem 8:} We have $$\mbox{End}_{\mathfrak{g}}(\mathbb{W}) = \mathfrak{C}.$$

{\it Proof.}  
Due to  Lemma 5, we know that $\mathfrak{C} \subseteq
\hbox{End}_{\mathfrak{g}}(\bigwedge^{\bullet}(\mathbb{V}^*)^{\mathbb{C}}\otimes \mathbb{S}).$ We prove the opposite inclusion. 
Obviously, we have $\mbox{End}_{\mathfrak{g}}(\mathbb{W})=
\mbox{End}_{\mathfrak{g}}(\bigwedge^{\bullet}{(\mathbb{V}^*)}^{\mathbb{C}}\otimes \mathbb{S}_+)\oplus
\mbox{End}_{\mathfrak{g}}(\bigwedge^{\bullet}{(\mathbb{V}^*)}^{\mathbb{C}}\otimes \mathbb{S}_-)\oplus
\mbox{Hom}_{\mathfrak{g}}(\bigwedge^{\bullet}{(\mathbb{V}^*)}^{\mathbb{C}}\otimes \mathbb{S}_+,\bigwedge^{\bullet}{(\mathbb{V}^*)}^{\mathbb{C}}\otimes \mathbb{S}_-)
\oplus \mbox{Hom}_{\mathfrak{g}}(\bigwedge^{\bullet}{(\mathbb{V}^*)}^{\mathbb{C}}\otimes \mathbb{S}_-,\bigwedge^{\bullet}{(\mathbb{V}^*)}^{\mathbb{C}}\otimes \mathbb{S}_+).$ Let us 
prove $\mbox{Hom}_{\mathfrak{g}}(\bigwedge^{\bullet}{(\mathbb{V}^*)}^{\mathbb{C}}\otimes \mathbb{S}_+,\bigwedge^{\bullet}{(\mathbb{V}^*)}^{\mathbb{C}}\otimes \mathbb{S}_-)\subseteq \mathfrak{C}$ only.
Due to the structure of the decomposition of $\bigwedge^{\bullet}(\mathbb{V}^{*})^{\mathbb{C}}\otimes \mathbb{S}$  (Theorem 3), we see that
$\hbox{Hom}_{\mathfrak{g}}(\bigwedge^{\bullet}{(\mathbb{V}^*)}^{\mathbb{C}}\otimes \mathbb{S}_+,
\bigwedge^{\bullet}{(\mathbb{V}^*)}^{\mathbb{C}}\otimes \mathbb{S}_-)\simeq \bigoplus_{(i,j), (k,m)\in \Xi}
\hbox{Hom}_{\mathfrak{g}}(\mathbb{E}_{ij}^+, \mathbb{E}_{km}^-).$   
Regarding the highest weights of the irreducible $\mathfrak{g}$-modules $\mathbb{E}_{ij}^{+}$ and $\mathbb{E}_{km}^{-}$ and using the 
 the   Schur lemma, we see that the following is true. If there exists a non trivial $\mathfrak{g}$-equivariant map between $\mathbb{E}_{ij}^+$ and $\mathbb{E}_{km}^{-},$ then $j=m$ and $k= 2r+i+1$ for some $r \in \mathbb{Z}.$
Suppose $r\geq 0.$ 
Due to Lemma 7,  Lemma 6 and the remark below this lemma,
we know that $(F^{+})^{1+2r}S_{ij}^{+}: \mathbb{E}_{ij}^+
\to \mathbb{E}_{i+2r+1,j}^-$  for each $r\in \mathbb{N}_0$ and $(i,j)\in \Xi.$
Due to  Lemma 7 again, we know that the mapping is nontrivial whenever $(i,j)\in \Xi_+$ and $(i+2r+1,j) \in \Xi.$
Thus $\mathbb{C}(F^{+})^{1+2r}S_{ij}^+=\mbox{Hom}_{\mathfrak{g}}(\mathbb{E}_{ij}^+,\mathbb{E}_{i+1+2r,j}^-)$ due to the   Schur lemma.
Similarly one can proceed for $r<0$ and $F^-.$ 
Using the fact that $S_{ij}^{+}\in \mathfrak{C}$ (Lemma 7), the inclusion in the considered special case follows. Tracing the remaining cases, 
we get $\mathfrak{C}= \mbox{End}_{\mathfrak{g}}(\mathbb{W})$.
$\Box$



{\bf Remark:} Let us remark that more conceptually, we could have defined the projections $S_{ij}^{\pm}$ onto the submodules $\mathbb{E}_{ij}^{\pm}$ using the Casimir operator of the super Lie algebra $\mathfrak{osp}(1|2).$

Now, we  introduce a representation of the super ortho-symplectic Lie algebra $\mathfrak{osp}(1|2)=:\mathfrak{g}'=\mathfrak{g}_0'\oplus \mathfrak{g}_1'$ on $\mathbb{W}$ 
 in order to decompose the space $\mathbb{W}$   as a $(\mathfrak{g} \times \mathfrak{g}')$-module.
The super Lie bracket of two $\mathbb{Z}_2$-homogeneous elements $u, v \in \mathfrak{g}'$ will be denoted by $[u,v]$ if and only if at least one of them is from the even part $\mathfrak{g}'_0.$ In the other cases, we will denote it by $\{u,v\}.$ Further there exists a basis $\{ h, e^+, e^-, f^+, f^-\}$ of $\mathfrak{g}',$ such that the set $\{e^+,h,e^-\}$ spans the even part $\mathfrak{g}'_0,$ the set $\{f^+,f^-\}$ spans the odd part $\mathfrak{g}'_1$ and
the only nonzero relations among the basis elements are
\begin{equation}
[h,e^{\pm}]= \pm e^{\pm} \qquad [e^+,e^-]=2h \label{he}
\end{equation}
\begin{equation}
[h,f^{\pm}]=\pm \frac{1}{2}f^{\pm} \qquad \{f^{+},f^-\}=\frac{1}{2}h \label{hf}
\end{equation}
\begin{equation}
[e^{\pm},f^{\mp}]=-f^{\pm} \qquad \{f^{\pm},f^{\pm}\}=\pm \frac{1}{2}e^{\pm}. \label{ef}
\end{equation}

For $r=0,\ldots, 2l,$ let us introduce the operators 
$$E^{\pm}:\bigwedge ^{r}(\mathbb{V}^*)^{\mathbb{C}}\otimes \mathbb{S} \to \bigwedge ^{r\pm 2}(\mathbb{V}^*)^{\mathbb{C}}\otimes \mathbb{S}, \,
E^{\pm}=\pm 2 \{F^{\pm},F^{\pm}\}$$ where $\{,\}$ denotes the anticommutator in the associative algebra $\mbox{End}(\mathbb{W}).$
Consider the following
mapping $\sigma: \mathfrak{osp}(1|2) \to  \hbox{End}(\bigwedge^{\bullet}(\mathbb{V}^*)^{\mathbb{C}}\otimes \mathbb{S})$ defined by
$$\sigma(e^{\pm}):=E^{\pm}, \, \mbox{ } \sigma(f^{\pm}):=F^{\pm} \, \mbox{ and } \sigma(h):=H$$
and extend it linearly to the whole $\mathfrak{g}'=\mathfrak{osp}(1|2).$

We prove that $\sigma$ is a super Lie algebra representation.
To do it, let us introduce a $\mathbb{Z}_2$-grading on the vector space $\mathbb{W}.$
We set $\mathbb{W}_0 :=(\bigoplus_{i = 0}^{l}\bigwedge^{2i}(\mathbb{V}^*)^{\mathbb{C}}) \otimes \mathbb{S}$ and $\mathbb{W}_1 :=(\bigoplus_{i=0}^{l-1}\bigwedge^{2i+1}(\mathbb{V}^*)^{\mathbb{C}}) \otimes \mathbb{S}.$ The vector space $\mbox{End}(\bigwedge^{\bullet}(\mathbb{V}^{*})^{\mathbb{C}} \otimes \mathbb{S})$ will be considered with the super Lie algebra structure inherited from the super vector space structure $\mathbb{W}=\mathbb{W}_0\oplus \mathbb{W}_1.$  

{\bf Theorem 9:} The mapping $$\sigma: \mathfrak{osp}(1|2)\to \hbox{End}(\bigwedge^{\bullet}(\mathbb{V}^*)^{\mathbb{C}}\otimes \mathbb{S})$$
is a super Lie algebra representation.

{\it Proof.} First, it is easy to see that $\sigma(\mathfrak{g}'_i) \subseteq \mbox{End}_i(\mathbb{W}),$ $i=0,1.$
Second, we shall verify that the operators $E^{\pm}, H, F^{\pm}$ satisfy the same relations as the ones for $e^{\pm}, h, f^{\pm}$ written in the rows (\ref{he}), (\ref{hf}), (\ref{ef}) above.
For $r=0,\ldots, 2l$ and $\alpha \otimes s \in \bigwedge^r (\mathbb{V}^*)^{\mathbb{C}} \otimes \mathbb{S},$ we have
\begin{eqnarray*}
[H,F^{+}](\alpha \otimes s)&=& H F^+ (\alpha \otimes s) - F^+ H(\alpha \otimes s)\\
&=& H (\frac{\imath}{2}\sum_{i=1}^{2l}\epsilon^i\wedge \alpha \otimes e_i .s) - F^+ \frac{1}{2}(r-l)(\alpha \otimes s)\\
&=& \sum_{i=1}^{2l}\left[\frac{1}{2}\frac{\imath}{2} (r+1-l)\epsilon^i \wedge \alpha \otimes e_i.s  - \frac{\imath}{2}\frac{1}{2}(r-l)\epsilon^i \wedge \alpha \otimes e_i.s\right]\\
&=&\frac{\imath}{4} \sum_{i=1}^{2l}\epsilon^i \wedge \alpha \otimes e_i.s = \frac{1}{2}F^+ (\alpha \otimes s).
\end{eqnarray*}
Thus we got the relation in the form of the $+$ version of the first equation written in the row $(\ref{hf})$ as required. Similarly, one can prove the $-$ version of the first equation written in that row. The second relation written in the row $(\ref{ef})$ and  the second relation in the row $(\ref{hf})$ follow easily from the definitions of $E^{\pm}$ and $H,$ respectively. The remaining relations, i.e., that ones in the row (\ref{he}) and the first relation in the row $(\ref{ef})$ can be proved using just the derived ones and by expanding the commutator and anticommutator of compositions of endomorphisms. We shall show explicitly how to prove the first relation in the row (\ref{ef}) only. Using   the definitions of the considered mappings, we may write 
$[E^+,F^-]=[2\{F^+,F^+\}, F^-]=4[F^+F^+, F^-]=4 (F^+F^+F^- - F^-F^+F^+) = 4 [F^+(-F^-F^+ + \frac{1}{2}H)-F^-F^+F^+]=
4 (F^-F^+F^- - \frac{1}{2}HF^+ + \frac{1}{2}F^+H- F^-F^+F^+)=2 [F^+,H] = - F^-.$
$\Box$

Summing up, we have the following

{\bf Corollary 10:} The super Lie algebra representation $$\sigma: \mathfrak{osp}(1|2) \to \mbox{End}(\bigwedge^{\bullet}{(\mathbb{V}}^*)^{\mathbb{C}}\otimes \mathbb{S})$$ maps the super Lie algebra
$\mathfrak{osp}(1|2)$ into the commutant algebra  $\mbox{End}_{\mathfrak{g}}(\mathbb{W})$ of $\mathfrak{g}$-invariants.

{\it Proof.} Follows from Lemma 5 and Theorem 9. $\Box$

\section{Howe duality for $Mp(\mathbb{V},\omega)$ acting on $\mathbb{W}$}
 
 We start by defining a certain family $\{\sigma_j\}_{j=0}^l$ of finite dimensional  irreducible representations of the super Lie ortho-symplectic algebra $\mathfrak{osp}(1|2).$
   For $j=0,\ldots, l,$ let $\mathbb{G}^j$ be a complex vector space of complex dimension $2l-2j+1$ and consider a basis $\{f_{i}\}_{i=j}^{2l-j}$ of $\mathbb{G}^j.$  
The super vector space structure on $\mathbb{G}^j$ is defined as follows. For $j=0,\ldots, l,$ we set  $(\mathbb{G}^j)_0:=\hbox{Span}_{\mathbb{C}}(\{f_i| i \in \{j,\ldots, 2l-j\}\cap 2\mathbb{N}_0  \})$ and complementarly, we define $(\mathbb{G}^j)_1:=\hbox{Span}_{\mathbb{C}}(\{f_i| i \in \{j,\ldots, 2l-j\} \cap (2\mathbb{N}_0 + 1)\}).$
For convenience, we set $f_{k}:=0$ for $k \in \mathbb{Z} \setminus \{j,\ldots, 2l-j\}.$ 
We will not  denote the dependence of the basis elements on the number $j$ explicitly.   As a short hand for each $(i,j)\in \Xi,$ we introduce the 
$\mathbb{Q}$-numbers
$$A(l,i,j):=\frac{(-1)^{i-j}+1}{16}(i-j) + \frac{(-1)^{i-j+1}+1}{16}(i + j - 2l - 1).$$

For $j=0,\ldots, l,$ we define
the mentioned representations 
$$\sigma_j: \mathfrak{osp}(1|2)\to \hbox{End}(\mathbb{G}^j)$$
\begin{eqnarray*}
\sigma_j(f^+)(f_i)&:=&A(l,i+1,j)f_{i+1},\, i = j,\ldots, 2l-j,\\
\sigma_j(f^-)(f_i)&:=&f_{i-1}, \, i=j,\ldots, 2l-j,\\
\sigma_j(h)&:=&2\{\sigma_j(f^+), \sigma_j(f^{-})\}, \\
\sigma_j(e^{\pm})&:=&\pm 2\{\sigma_j(f^{\pm}),\sigma_j(f^{\pm})\}.
\end{eqnarray*}
  
We have the following  

 {\bf Lemma 11:} For $j=0,\ldots, l,$ the mapping $\sigma_j: \mathfrak{osp}(1|2)\to \hbox{End}(\mathbb{G}^j)$ is an irreducible representation of the super Lie algebra $\mathfrak{osp}(1|2).$ 

{\it Proof.}
First, we prove that for $j=0,\ldots, l,$ the mapping  $\sigma_j$ is a   representation of the super Lie algebra $\mathfrak{osp}(1|2).$
It is easy to see that  whereas the even part of $\mathfrak{g}'$ acts by transforming the even part of $\mathbb{G}^j$ into itself and the odd part into itself  as well,  the odd part of $\mathfrak{g}'$ acts by interchanging the mentioned two parts of $\mathbb{G}^j.$ 

Now, we should check weather the relations (\ref{he}), (\ref{hf}) and (\ref{ef}) are preserved by the mapping $\sigma_j$ for $j=0,\ldots, l.$
The  last relation in the row (\ref{ef}) and the last relation in (\ref{hf}) are satisfied due the definitions of $\sigma_j(e^{\pm})$ and $\sigma_j(h),$ respectively.
Let us prove the $+$ version of the first relation written in the row (\ref{hf}).
For $i=j,\ldots, 2l-j,$ the left hand side of the mentioned relation reads.
\begin{multline*}
 [\sigma_j(h)\sigma_j(f^+)-\sigma_j(f^+)\sigma_j(h)]f_i\\
 \begin{aligned}
= & 2[(\sigma_j(f^+)\sigma_j(f^-)+ \sigma_j(f^-)\sigma_j(f^+))\sigma_j(f^+)\\
  & - \sigma_j(f^+)(\sigma_j(f^+)\sigma_j(f^-)+\sigma_j(f^-)\sigma_j(f^+)])f_i\\
= & 2[\sigma_j(f^-) \sigma_j(f^+)\sigma_j(f^+) -\sigma_j(f^+) \sigma_j(f^+) \sigma_j(f^-)]f_i\\
= & 2[A(l,i+1,j)\sigma_j(f^-)\sigma_j(f^+)f_{i+1}-\sigma_j(f^+)\sigma_j(f^+) f_{i-1}]\\
= & 2[A(l,i+2,j)A(l,i+1,j) - A(l,i,j)A(l,i+1,j)]f_{i+1}\\
= & 2A(l,i+1,j)(A(l,i,j) - A(l,i+1,j))f_{i+1}\\
= & \frac{1}{2}A(l,i+1,j)f_{i+1}=\frac{1}{2}\sigma_j(f^+)f_{i}.
\end{aligned}
\end{multline*}
The $-$ version of this relation can be proved in a similar way. To verify the relations written in the row (\ref{he}) and the first relation in the row (\ref{ef}), it is sufficient to use the already proved relations and expand the definition of commutator and anticommutator of compositions of endomorphisms.

  Now, we prove that for $j=0,\ldots, l,$ the representation $\sigma_j$ is irreducible.
Let us suppose, there exists a non-trivial proper invariant subspace $X \subseteq \mathbb{G}^j$ of dimension $0<k<2l-2j+1.$
Let us chose a basis $\{v_j,\ldots, v_{j+k-1}\}$ of $X.$  Let $Y:=\hbox{Span}_{\mathbb{C}}(\{f_j,\ldots, f_{j+k-1}\})$ and define a vector space isomorphism $T: X \to Y$ by the prescription
$T(v_i):=f_i,$ $i=j,\ldots, j+k-1$ and extend it first linearly and then to an automorphism $\tilde{T}$ of $\mathbb{G}^j.$
Now, for each $j=0,\ldots,l$, consider a representation $\tilde{\sigma}_j: \mathfrak{osp}(1|2) \to \hbox{End}(\mathbb{G}^j)$ given by $\tilde{\sigma}_j(X):=\tilde{T}^{-1}\sigma_j(X)\tilde{T}$ for each $X\in \mathfrak{osp}(1|2).$ This representation is clearly equivalent to $\sigma_j.$
 Using the definition of $\sigma_j,$  we get $\tilde{\sigma_j}(f^+)v_{j+k-1}=\tilde{T}^{-1}\sigma_j(f^+)\tilde{T} v_{j+k -1}
=\tilde{T}^{-1}\sigma_j(f^+)f_{j+k-1} = \tilde{T}^{-1}f_{j-k} \notin X,$  i.e.,  $X$ is not invariant.
Thus we see, that $\tilde{\sigma}_j$ is irreducible for $j = 0 ,\ldots, l.$ Because of the equivalence
of $\tilde{\sigma}_j$ and $\sigma_j,$ the lemma follows.
$\Box$

Now, we prove a technical but quite important

{\bf Lemma 12:} For each $k \in \mathbb{N}_0$ and $i=0,\ldots, 2l,$ we have 
\begin{eqnarray*}
&&(F^-)^k F^+=\\
&&(-1)^k F^+(F^-)^k + \left[\frac{(-1)^k+1}{16}k+\frac{(-1)^{k+1}+1}{16}(2i-2l-k+1)\right](F^{-})^{k-1},
\end{eqnarray*}
when acting on $\bigwedge^i(\mathbb{V}^*)^{\mathbb{C}}\otimes \mathbb{S}.$

{\it Proof.} We will not write explicitly that we are considering the action of the considered operators on the space $\bigwedge^{i}(\mathbb{V}^*)^{\mathbb{C}}\otimes \mathbb{S}$ and proceed by induction.
\begin{itemize}
\item[I.] For $k=0$ the lemma holds obviously.
\item[II.] 
\begin{itemize}
\item[a.]We suppose that the lemma holds for an even integer $k\in \mathbb{N}_0.$
We have 
\begin{eqnarray*}
(F^-)^{k+1}F^+&=& F^-(F^-)^kF^+\\
              &=& F^-[F^+(F^-)^k+\frac{k}{16}((-1)^k+1)(F^-)^{k-1}]\\
              &=& -(F^+)(F^-)^{k+1}+\frac{1}{2}H(F^-)^k+2\frac{k}{16}(F^{-})^k\\
              &=& -F^+(F^-)^{k+1}+\frac{1}{4}(i-k-l)(F^{-})^{k}+\frac{k}{8}(F^{-})^k\\
              & & -F^+(F^-)^{k+1}+\frac{2}{16}(2i-(k+1)+1-2l)(F^{-})^k,
\end{eqnarray*}
where we have  used the induction hypothesis, commutation relation $2\{F^{+}, F^{-}\}=H$
and  Lemma 4 on the value of $H.$ The last written expression coincides with that one of statement of the lemma for $k+1$ odd.
\item[b.] Now, suppose $k$ is odd.
We have
\begin{eqnarray*}
&&(F^-)^{k+1}F^{+}=F^-(F^-)^{k}F^{+}\\
&=&F^-[-F^+(F^-)^{k}+\frac{(-1)^{k+1}+1}{16}(2i-2l-k+1)(F^{-})^k]\\
&=&+F^+(F^-)^{k+1}-\frac{1}{2}H(F^{-})^k+\frac{1}{8}(2i-2l-k+1)(F^{-})^k\\  
 &=&F^+(F^-)^{k+1}-\frac{1}{8}(2i-2k-2l)(F^-)^k+\\
&&+\frac{1}{8}(2i-2l-k+1)(F^{-})^k =F^+(F^-)^{k+1}-\frac{2}{16}(k+1)(F^{-})^k,
\end{eqnarray*}               
where  we have again used the same tools as in the previous item.
\end{itemize}
\end{itemize}
$\Box$

Now let us define a family $\{\rho_j^{\pm} | j=0,\ldots, l\}$ of representations of the symplectic Lie algeba $\mathfrak{g}$ on $\mathbb{E}_{jj}^{\pm}$  by  the prescription
$$\rho_j^{\pm}: \mathfrak{g} \to \mbox{End}(\mathbb{E}_{jj}^{\pm}) \mbox{  } \, ,
\rho_j(X)v:=\rho(X)v$$ for $X\in \mathfrak{g}$ and $v \in \mathbb{E}_{jj}^{\pm}.$
  Let us introduce the following mapping $\mbox{Sgn}: \{\pm\} \times \mathbb{N}_0 \to \{\pm\}$ given by the prescription $\mbox{Sgn}(\pm, 2k):=\pm$ and $\mbox{Sgn}(\pm,2k+1)=\mp$ for each $k\in \mathbb{Z}.$
Now for $(i,j) \in \Xi,$ let us define a mapping
$\psi_{ij}^{\pm}: \mathbb{E}_{ij}^{\pm} \to \mathbb{E}_{jj}^{\mbox{Sgn}(\pm,i-j)} \otimes \mathbb{G}^j$
by the formula 
$$\psi_{ij}^{\pm} v: = (F^{-})^{i-j}v\otimes f_i$$ for
an element $v \in \mathbb{E}_{ij}^{\pm}.$ Further we set $\psi:=\bigoplus_{(i,j)\in \Xi}(\psi_{ij}^+ \oplus \psi_{ij}^-).$
Let us  consider $\bigwedge^{\bullet}(\mathbb{V}^*)^{\mathbb{C}}\otimes \mathbb{S}$ with the action $L = \rho \otimes \sigma$ of $\mathfrak{g}\times \mathfrak{g}'$ and  the space $\bigoplus_{j=0}^l(\mathbb{E}_{jj}^{+}\oplus \mathbb{E}_{jj}^-) \otimes \mathbb{G}^j$ with the action $R =  \bigoplus_{j=0}^l(\rho_j^+ \oplus \rho_j^-) \otimes \sigma_j$ of $\mathfrak{g} \times \mathfrak{g}'.$  

In the next theorem, the aforementioned Howe duality is stated. 

{\bf Theorem 13:} We have the following $(\mathfrak{sp}(\mathbb{V}^{\mathbb{C}},\omega) \times \mathfrak{osp}(1|2))$-module isomorphism
$$\bigwedge^{\bullet}(\mathbb{V}^*)^{\mathbb{C}}\otimes \mathbb{S} \simeq \bigoplus_{j=0}^{l}
(\mathbb{E}_{jj}^+ \oplus \mathbb{E}_{jj}^-)\otimes \mathbb{G}^j.$$
{\it Proof.} 
Due to Theorem 3 we know that
$\bigwedge^{\bullet}(\mathbb{V}^*)^{\mathbb{C}}\otimes \mathbb{S} = \bigoplus_{(i,j)\in \Xi}(\mathbb{E}_{ij}^+\oplus \mathbb{E}_{ij}^-).$ 
It is evident, that $\psi$ is a vector space isomorphism.
We prove that for each $(i,j)\in \Xi,$ the mapping $\psi_{ij}^{\pm}: \mathbb{E}_{ij}^{\pm} \to \mathbb{E}_{jj}^{\mbox{Sgn}(\pm,i-j)}\otimes \mathbb{G}^j$ is $(\mathfrak{g}\times \mathfrak{g}')$-equivariant. The $\mathfrak{g}$-equivariance follows easily because $F^-$ in the definition of $\psi_{ij}^{\pm}$ commutes with the representation $\rho$ of $\mathfrak{g}$, see Lemma 5. We shall check the $\mathfrak{g}'$-equivariance.
For each $(i,j)\in \Xi$ and $v \in \mathbb{E}_{ij}^{\pm},$ consider 
$\psi_{ij}^{\pm}\sigma(f^-)v=\psi_{ij}^{\pm}F^- v=(F^-)^{i-1-j}F^-v \otimes f_{i-1}=(F^-)^{i-j} v \otimes f_{i-1}.$
On the other hand, we have
$\sigma_j(f^-)(\psi_{ij}^{\pm}v)=F^-((F^{-})^{i-j}v\otimes f_i)=(F^{-})^{i-j}v\otimes f_{i-1}.$ 
Now, we check the $\mathfrak{g}'$-equivariance in the case of $F^+.$
We shall use Lemma 12 to compute $\psi_{ij}^{\pm}\sigma_j(f^+)v=\psi_{ij}^{\pm}F^+v=(F^-)^{i+1-j}F^+v\otimes f_{i+1}=[(-1)^{i+1-j}F^+(F^-)^{i+1-j}v+A(l,i+1,j)(F^-)^{i-j}v]\otimes f_{i+1}=
A(l,i+1,j)(F^{-})^{i-j}v \otimes f_{i+1},$ where we have used the fact that $(F^-)^{i+1-j}v=0$ because $v\in \mathbb{E}_{ij}^{\pm}.$
Let us compute $\sigma(f^+)\psi_{ij}^{\pm}v=F^+\psi_{ij}^{\pm}v=(F^-)^{i-j}v\otimes A(l,i+1,j)f_{i+1}.$ Thus the equivariance with respect to $F^+$ is proved.
Because the operators $H,$ $E^+$ and $E^-$ are linear combinations of compositions of the operators $F^+$ and $F^-,$ the $\mathfrak{g}'$-equivariance of $\psi_{ij}^{\pm}$ follows. $\Box$


{\bf Remark:} One may call this Howe correspondence 2-folded or $2:1$ because the pairing of $\mathbb{E}_{jj}^+$ and $\mathbb{E}_{jj}^{-}$ with $\mathbb{G}^j$ is not one to one. Namely, to $\mathbb{G}^j$ the module $\mathbb{E}_{jj}^+$ as well as the module $\mathbb{E}_{jj}^-$ are paired in the corresponding tensor product.
That is basically  due to the fact that we are decomposing with respect to
$\mathfrak{g} \times \mathfrak{g}'$ only not concerning the whole commutant of $\mathfrak{g}$ acting on $\mathbb{W}$.
 
{\bf Remark:} Due to the fact that the category of admissible finite length Harish-Chandra modules is a full subcategory of the category of $\mathcal{U}(\mathfrak{g})$-modules and some 
basic  properties of the minimal globalization functors, the  results of the paper have their appropriate minimal globalization counterparts.

\end{document}